\documentclass[a4paper,10pt]{amsart}
\usepackage{epsf}  
\usepackage{amsfonts, amssymb, amsthm, amsmath}  
\usepackage{hyperref}
\usepackage{tikz}
\usetikzlibrary{matrix,arrows,decorations.pathmorphing}
\usepackage{tikz-cd}
\usepackage{graphicx}
\usepackage{psfrag}
\DeclareMathOperator{\Aut}{Aut}

\newcommand{\rh}{{}^{r}H}
\newcommand{\ex}{\bold}

\providecommand{\C}[2]{\ensuremath {C^{#1,\underline{#2}}}}

\newcommand{\dbar}{\bar{\partial}}

\providecommand{\lrb}[1]{\ensuremath{\left(#1\right)}}
\providecommand{\abs}[1]{\left\lvert #1\right\rvert}

\newcommand{\R}{\bold{\mathcal R}}

\thanks{Supported by the Australian Research Council grant DP140100296. I also wish to thank Aleksey Zinger for convincing me to write this case of the gluing formula separately, and for comments on an earlier version.}

\author{Brett Parker   }
\email{brettdparker@gmail.com}
\title{Gluing formula for Gromov--Witten invariants in a triple product}

\begin{document}

\begin{abstract} We present a gluing formula for Gromov--Witten invariants in the case of a triple product. This gluing formula is a simple case of a general tropical gluing formula proved and stated using exploded manifolds. We present this simple case because it is relatively easy to explain without any knowledge of exploded manifolds or log schemes. 
\end{abstract}

\maketitle

\tableofcontents

\section{Introduction}
The symplectic-sum formula, \cite{Li,IP,ruan}, is an important tool for calculating Gromov--Witten invariants. In this paper, we explain its simplest generalization: a gluing formula for a triple product. This simplest generalization has surprising complications not present in earlier gluing formulae --- tropical curves play a staring role, and we can not always use ordinary cohomology. Once these complications are  mastered, however, there are no further surprises in the general case of a gluing formula for Gromov--Witten invariants using an arbitrary normal-crossing or log-smooth degeneration.

Consider a complex\footnote{We work in the complex category for expositional convenience. The same gluing formula holds in the  symplectic analogue of this setting when almost complex structures are chosen suitably compatible with divisors.} family of compact K\"ahler manifolds over a complex disk, smooth away from the origin, over which there are  singularities  at worst in the form $z_{1}z_{2}z_{3}=t$. Assume furthermore that the fiber over the origin is the union of $3$ connected K\"ahler manifolds $M_{1},M_{2},M_{3}$, and that all other strata --- $M_{i}\cap M_{j}\setminus (M_{1}\cap M_{2}\cap M_{3})$ and $M_{1}\cap M_{2}\cap M_{3}$ --- are also connected. Over $M_{i}\cap M_{j}$, there is a natural $\mathbb C^{*}$--bundle, $L_{ij}$, equal to the normal bundle of $M_{i}\cap M_{j}\subset M_{i}$ minus its zero section.  

This paper describes how to compute the Gromov--Witten invariants of a generic fiber $M$ of this family in terms of some relative Gromov--Witten invariants of $M_{i}$, $L_{ij}$, and $L_{12}\oplus L_{13}\rvert_{M_{1}\cap M_{2}\cap M_{3}}$.\footnote{By $L_{12}\oplus L_{13}\rvert_{M_{1}\cap M_{2}\cap M_{3}}$ we mean the fiber product of $L_{12}\rvert_{M_{1}\cap M_{2}\cap M_{3}}$ with $L_{13}\rvert_{M_{1}\cap M_{2}\cap M_{3}}$ over $M_{1}\cap M_{2}\cap M_{3}$.} To avoid discussion of the relationship between homology of a generic fiber with the homology of our pieces $M_{i}$, we shall only consider point constraints, and we shall not keep track of the full homology class represented by curves; instead we  keep track of the energy of curves, defined by the integral of the K\"ahler form. We begin by describing the gluing formula in the genus $0$ case, then explain what complications arise for higher genus and families with more complicated normal crossing singularities. We follow this with sections \ref{CP2} and \ref{unbalanced}, which contain simple examples where consideration of tropical curves is essential.

Aspects of this gluing formula will surprise some readers familiar with  symplectic field theory or the symplectic sum formula for Gromov--Witten invariants. Why should tropical curves be relevant? Why should the contributions of relative invariants from $L_{ij}$ and $L_{12}\oplus L_{13}\rvert_{M_{1}\cap M_{2}\cap M_{3}}$ be so important? Why should curves with such ungeneric-looking domains contribute as in the example at the end of section \ref{CP2}? Why can't we just use usual cohomology? These oddities are unavoidable consequences of generalizing the symplectic-sum formula to a triple product, but become natural when we view the problem from the perspective of exploded manifolds or log geometry --- we do not use exploded manifolds in this paper, but outline how they help in Section \ref{exploded}.

\section{Tropical curves in the dual intersection complex} 

\

We shall need the dual intersection complex $\triangle$ of $\{M_{1},M_{2},M_{3}\}$, represented as the closed triangle in $\mathbb R^{2}$ with vertices $(0,0)$ corresponding to $M_{1}$, $(1,0)$ corresponding to $M_{2}$,  and $(0,1)$ corresponding to $M_{3}$.\footnote{The reader may object that the dual intersection complex does not look symmetric. In fact, we shall be using the integral-affine structure of $\mathbb R^{2}$, consisting of the lattice of integer vectors inside its tangent space. Any translation, and any invertible integral linear  transformation is a symmetry of $\mathbb R^{2}$ with this structure, so the $\triangle$ with this integral-affine structure is maximally symmetric.} A picture of $\triangle$ is below.

\includegraphics{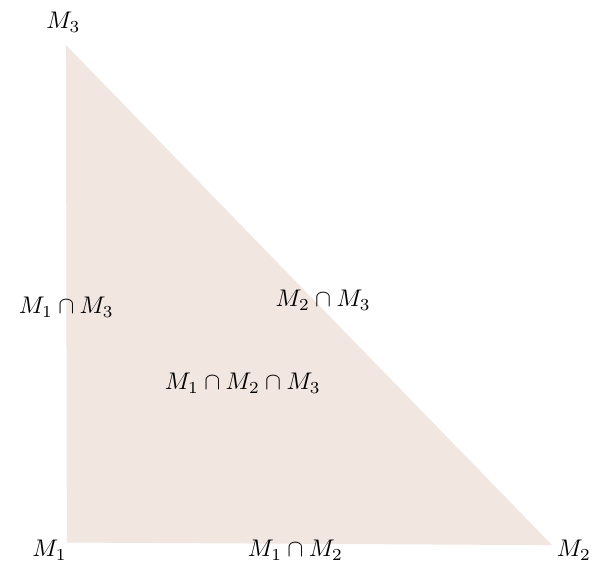}\\

\noindent  Our gluing formula consists of a sum over tropical curves in $\triangle$. By a tropical curve $\gamma$, we mean a continuous map of a complete metric graph into $\triangle$ so that this map restricted to each edge has integral derivative.  
  
 \includegraphics{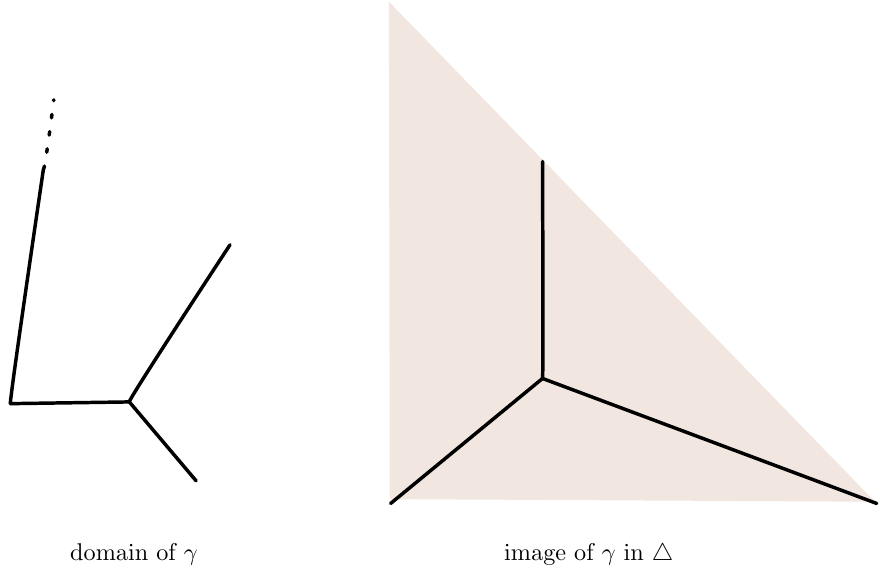}

\noindent More specifically, the domain of a tropical curve $\gamma$ consists of the following:
\begin{itemize}
\item a finite set of edges, isometric\footnote{In dimension $1$, an integral-affine structure is the same thing as a metric --- the integral vectors are the vectors with integer length.} to closed intervals in $\mathbb R$ with positive length, 
\item  a finite set of vertices,
\item a map from the boundary of each edge to the set of vertices; this map is used to glue edges to vertices to obtain a topological space --- we shall generally make the assumption that this topological space is connected, and say our tropical curve is connected. 
\end{itemize}
A tropical curve $\gamma$ in the dual intersection complex $\triangle$ is then a continuous map from the domain of $\gamma$ to $\triangle$,  integral-affine  when restricted to each edge in the sense that the derivative on each edge sends vectors of integer length to integral vectors in $\mathbb Z^{2}\subset \mathbb R^{2}$. In Section \ref{relative GW section}, we  see that the derivative $(a,b)\in \mathbb Z^{2}$ along an (oriented) edge is data describing the degree of contact of a holomorphic  curve with certain divisors. 

\section{Target stacks $\bar X_{e}$ for relative invariants.}

  For each edge $e$ of a tropical curve, there is an associated stack\footnote{The reader need not be familiar with stacks, as $\bar X_{e}$ will always be the quotient of a smooth manifold by a trivial group action, and we shall say explicitly what effect this has on cohomology.} $X_{e}$ with compactification $\bar X_{e}$ and the following significance:  the moduli space of holomorphic curves we want is a fiber product\footnote{Actually, this fiber product takes place in the category of exploded manifolds, which is why we need a different cohomology theory in the general case.} of relative moduli spaces over $\bar X_{e}$ for all internal edges $e$. The $0$-genus gluing formula will involve cohomology classes in $H^{*}(\bar X_{e},\mathbb R)$, and the higher genus gluing formula will involve cohomology classes in $\rh^{*}(\prod_{e}\bar X_{e},\mathbb R)$, where $\rh$ is refined cohomology, discussed further in section \ref{refined cohomology}.

 \begin{itemize}
 \item If $e$ is contained in the interior of the dual intersection complex $\triangle$ and has derivative $(a,b)$ then $X_{e}$ is the quotient of $L_{12}\oplus L_{13}\rvert_{M_{1}\cap M_{2}\cap M_{3}}$ by the $\mathbb C^{*}$--action with weight $(a,b)$. If $(a,b)$ is $k\geq 1$ times a primitive integral vector, then this  $X_{e}$ the quotient of  a $\mathbb C^{*}$--bundle over $M_{1}\cap M_{2}\cap M_{3}$ by a trivial $\mathbb Z_{k}$--action, and $\bar X_{e}$ is the quotient of the associated $\mathbb CP^{1}$--bundle by the trivial $\mathbb Z_{k}$--action.  If $(a,b)=(0,0)$, then we can pick a fiberwise-toric compactification  of $L_{12}\oplus L_{13}\rvert_{M_{1}\cap M_{2}\cap M_{3}}$, then define $\bar X_{e}$ to be the quotient of this compactification by the trivial $\mathbb C^{*}$--action.
 
 \item If $e$ is contained in the boundary of the dual intersection complex corresponding to $M_{i}\cap M_{j}$, and has derivative $k\geq 1$ times a primitive integral vector, then $\bar X_{e}$ is the quotient  of $M_{i}\cap M_{j}$ by the trivial $\mathbb Z_{k}$--action. $X_{e}$ is obtained from this by removing its intersection with the remaining $M_{k}$. If such an edge  has derivative  $0$, then $\bar X_{e}$ is the quotient of the $\mathbb CP^{1}$--bundle containing $L_{ij}$ by the trivial $\mathbb C^{*}$--action\footnote{Although $\bar X_{e}$ will be different if $L_{ji}$ is used instead of $L_{ij}$, this will not affect anything. }, and $X_{e}$ is obtained by removing the intersection with the remaining $M_{k}$, and also removing the $0$ and $\infty$ sections.
 \item If $e$ is contained in the corner associated to $M_{i}$, then $\bar X_{e}$ is the quotient of $M_{i}$ by the trivial $\mathbb C^{*}$--action, and $X_{e}$ is obtained by removing the intersection with the remaining $M_{j}$ and $M_{k}$. 
 \end{itemize}
 
 In each case, $\bar X_{e}$ is the quotient of a smooth manifold by a trivial group action.  When this group is finite $H^{*}(\bar X_{e},\mathbb R)$ is just the cohomology of this manifold; otherwise,  this group is $\mathbb C^{*}$, and  $H^{*}(\bar X_{e},\mathbb R)$ is the tensor product  of the cohomology of this manifold  with the $\mathbb C^{*}$--equivariant cohomology of a point, equal to $\mathbb R[u]$. Recall that the stack $\cdot/\mathbb C^{*}$ classifies line bundles,\footnote{See \cite{BehrendCOS} for an introduction to the cohomology of stacks.}  and pulling back $u$ over a map $Y\longrightarrow \cdot/\mathbb C^{*}$ gives the first Chern class of the corresponding line bundle over $Y$.

\section{Genus $0$ case}  

Use notation 
\[\langle\tau_{d_{1}}(p_{1})\dotsb\tau_{d_{n}}(p_{n})\rangle_{g,E}\]
to indicate the Gromov--Witten invariant of $M$ defined by integrating $\prod_{i}\psi_{i}^{d_{i}}$ over the moduli space  of  curves with genus $g$, energy $E$, and $n$ labeled marked points sent to $n$ chosen points in $M$, where $\psi_{i}$ indicates the first Chern class of the tautological line bundle associated to the $i$th marked point.  It shall be useful to combine the zero-genus Gromov--Witten invariants for different $E$ into the following series involving a dummy variable $q$.
\[\langle\tau_{d_{1}}(p_{1})\dotsb\tau_{d_{n}}(p_{n})\rangle_{0}:=\sum_{E\geq 0}\langle\tau_{d_{1}}(p_{1})\dotsb\tau_{d_{n}}(p_{n})\rangle_{0,E}q^{E}\]
Use the notation $\bold R$ to indicate the Novikov ring of series $\sum_{\lambda\in[0,\infty)} a_{\lambda}q^{\lambda}$ where $a_{\lambda}\in \mathbb R$ is nonzero for only finitely many $\lambda$ less than any chosen $N$.

Although this Gromov--Witten invariant does not depend on the particular choice of $n$ points in $M$, our gluing formula shall depend on a choice of $n$ points $p_{1},\dotsc ,p_{n}$ in the dual intersection complex $\triangle$ --- different choices of points give different decompositions of the same invariants. After these points are chosen,  our Gromov--Witten invariant decomposes into terms corresponding to tropical curves $\gamma$ whose domains are connected graphs with genus $0$, and  exactly $n$ semi-infinite ends.\footnote{We call an edge of $\gamma$ isometric to $[0,\infty)$ an end, or external edge of $\gamma$.} Furthermore,  these ends of $\gamma$ are labeled $1,\dotsc,n$ and sent to $p_{1},\dotsc,p_{n}$ respectively. This decomposition may be written as follows. 
\[\langle\tau_{d_{1}}(p_{1})\dotsb\tau_{d_{n}}(p_{n})\rangle_{0}=\sum_{\gamma}\frac 1{\abs{\Aut \gamma}}\langle\tau_{d_{1}}(p_{1})\dotsb\tau_{d_{n}}(p_{n})\rangle_{0,\gamma}\]
In the above, $\Aut \gamma$ indicates the group of automorphisms of $\gamma$ that preserve the labeling  of the ends. Only a finite set of tropical curves will have a nonzero contribution to the the coefficient of $q^{E}$ in the above sum --- constraints on which tropical curves contribute come from the relationship between the derivatives of edges of $\gamma$ and the contact data of relative invariants that contribute to the gluing formula, discussed in section \ref{relative GW section}; moreover,    the contribution of $\gamma$ will always be $0$ unless $\gamma$ is rigid in the sense that it can't be deformed within the space of tropical curves of the same combinatorial type. In particular, for a tropical curve to be rigid,  no internal edges can have derivative $0$, but all ends of $\gamma$ must have $0$ derivative in order to fit inside $\triangle$.
 
 Each $\langle\tau_{d_{1}}(p_{1})\dotsb\tau_{d_{n}}(p_{n})\rangle_{0,\gamma}$ is a kind of Gromov--Witten invariant computed from relative invariants by a gluing formula below.

  Choose an orientation of the edges of $\gamma$ so that each external edge is oriented inwards, and each vertex has at most one edge exiting it.   To each vertex $v$ with an edge $e_{v}$ exiting it, we can associate a relative Gromov--Witten invariant valued in $H^{*}(\bar X_{e_{v}},\bold R)$. 
  \[\eta_{0,\gamma_{v}}:\bigotimes_{\text{incoming }e}H^{*}(\bar X_{e},\bold R)\longrightarrow H^{*}(\bar X_{e_{v}},\bold R)\] The above tensor product is over all edges $e$ entering $v$. When $v_{0}$ is the vertex with no edges exiting it, we can associate a $\bold R$--valued invariant. 
 \[\eta_{0,\gamma_{v_{0}}}:\bigotimes_{\text{incoming }e}H^{*}(\bar X_{e},\bold R)\longrightarrow \bold R\]
We discuss these relative invariants further in section \ref{relative GW section}. Assign a homology class $\theta_{e}\in H^{*}(\bar X_{e},\bold R)$ to each edge $e$ as follows. 
\begin{enumerate}
\item If $e$ is the $i$th external edge, $\theta_{e}:=\tau_{d_{i}}(p_{i})\in H^{*}(\bar X_{e},\bold R)$ is $(-u)^{d_{i}}$ times the Poincare dual of a point, where $u$ is the generator of the $\mathbb C^{*}$-equivariant cohomology of a point.
\item If $e$ exits the vertex $v$, and $\theta_{e'}$ has already been defined for all edges $e'$ entering $v$,  then $\theta_{e}=\eta_{0,\gamma_{v}}(\otimes \theta_{e'})$.
\end{enumerate}
Once $\theta_{e}$ has been defined for all edges, the formula for our Gromov--Witten invariant is as follows:
\[\langle\tau_{d_{1}}(p_{1})\dotsb\tau_{d_{n}}(p_{n})\rangle_{0,\gamma}=\eta_{0,\gamma_{v_{0}}}\lrb{\bigotimes_{\text{incoming }e}\theta_{e}}\]
For some examples, see sections \ref{CP2} and \ref{unbalanced}.

The use of the Novikov ring $\bold R$ is just a trick to keep track of the energy of curves. The relative invariants decompose as $\eta_{0,\gamma_{v}}=\sum_{E}\eta_{0,E,\gamma_{v}}q^{E}$ where $\eta_{0,E,\gamma_{v}}$ uses $\mathbb R$ in place of $\bold R$, and involves curves with energy $E$. The coefficient of $q^{E}$ provided by the gluing formula then involves curves with total energy  $E$.

\section{Higher genus}

Use the following notation for the arbitrary genus Gromov--Witten invariant.
\[\langle\tau_{d_{1}}(p_{1})\dotsb\tau_{d_{n}}(p_{n})\rangle:=\sum_{g, E}\langle\tau_{d_{1}}(p_{1})\dotsb\tau_{d_{n}}(p_{n})\rangle_{g,E}q^{E}x^{2g-2+n}\]
This takes values in the ring $\R$ of formal Laurent series in $x$ with coefficients in $\bold  R$.  After choosing points $p_{1},\dotsc,p_{n}$ in the dual intersection complex, this Gromov--Witten invariant decomposes as a sum over connected rigid tropical curves $\gamma$ with $n$ labeled ends sent to $p_{1},\dotsc,p_{n}$ respectively. 
\[\langle\tau_{d_{1}}(p_{1})\dotsb\tau_{d_{n}}(p_{n})\rangle=\sum_{\gamma}\frac 1{\abs{\Aut \gamma}}\langle\tau_{d_{1}}(p_{1})\dotsb\tau_{d_{n}}(p_{n})\rangle_{\gamma}\]

This time, we consider ends $e$ attached to a vertex $v$ to be incoming, and each choice of an orientation $e'$ of an internal edge leaving $v$ to be an outgoing edge of $\gamma_{v}$ --- if an edge of $\gamma$ has two ends on $v$, there are two corresponding outgoing edges of $\gamma_{v}$.  Then we can associate   a relative Gromov--Witten invariant $\eta_{\gamma_{v}}$ to each vertex.
\[\eta_{\gamma_{v}}:\rh^{*}\lrb{\prod\bar X_{e},\R}\longrightarrow\rh^{*}\lrb{\prod\bar X_{e'},\R}\]
 The $\rh$ above indicates refined cohomology, which we shall describe briefly in section \ref{refined cohomology}. For now, the reader can treat it like usual cohomology, except $\rh^{*}(\prod_{e}\bar X_{e})$ is in general larger than $\otimes_{e}\rh^{*}(\bar X_{e})$.  By pushing forward the relevant $\tau_{d_{i}}(p_{i})$ classes using $\eta_{\gamma_{v}}$, and pulling back the result to $\prod_{e}\bar X_{e}$, we may define a class $\theta_{v}\in\rh^{*}(\prod_{e}\bar X_{e},\R)$.
 \[\theta_{v}:=\Delta^{*}\eta_{\gamma_{v}}\lrb{\prod_{i}\tau_{d_{i}}(p_{i})}\subset \rh^{*}(\prod_{e}\bar X_{e},\R)\] 
In the above, the product of cohomology classes $\tau_{d_{i}}(p_{i})$ is over all $i$ labeling the external edges attached to $v$,  the product of  $\bar X_{e}$ is over all internal edges $e$ of $\gamma$, and $\Delta$ indicates the diagonal map from $\prod_{e}\bar X_{e}$ to $\prod_{e}\bar X_{e}^{2}$ followed by projection to the $\prod_{e'}\bar X_{e'}$ found in the target of $\eta_{\gamma_{v}}$. We may then write the gluing formula as follows:
\[\langle\tau_{d_{1}}(p_{1})\dotsb\tau_{d_{n}}(p_{n})\rangle_{\gamma}=\int_{\prod_{e}\bar X_{e}}\prod_{v}\theta_{v}\]
Note that $\prod_{e}\bar X_{e}$ is a stack obtained by taking the quotient of some manifold $\bar X$ by a trivial $G$--action. Integration over $\prod_{e} \bar X_{e}$ corresponds to integration over $\bar X$ divided by $\abs G$. In particular, the above integral vanishes if any of the internal edges of $\gamma$ have zero derivative.   

Again, the use of the fancy ring $\R$ is simply a trick to keep track of energy and Euler characteristic. The coefficient of $x^{a}$  involves curves with Euler characteristic $-a$ (when the marked points are regarded as punctures). The effect of our use of $\R$ is that the genus $g$ contributions of the gluing formula involve curves at vertices of $\gamma$ with genus adding up to $g$ minus the genus of $\gamma$.

\

The generalization of this gluing formula to an arbitrary simple normal crossings degeneration is straightforward. The only difference is that the dual intersection complex in which we count curves is more complicated. The same gluing formula also works for  Gromov--Witten invariants relative simple normal crossing divisors --- the dual intersection complex of a manifold with normal crossing divisors is a complex with a copy of $[0,\infty)^{k}$ for every $k$-fold intersection of divisors.

 The $0$-genus gluing formula also generalizes, but refined cohomology must be used instead of ordinary cohomology whenever the dimension of the dual intersection complex is greater than $2$. We shall discuss this further in Section \ref{refined cohomology}. 

We may also use a similar gluing formula to calculate Gromov--Witten invariants with insertions of different cohomology classes, however at this point we need to be more careful about how cohomology classes in $M$ can correspond to cohomology classes in the singular fiber. This is best done using exploded manifolds, so I will not describe it here.

\section{Relative Gromov--Witten invariants $\eta_{\gamma_{v}}$}\label{relative GW section}
We will not give the full definition of the relative Gromov--Witten invariants $\eta_{\gamma_{v}}$; see \cite{evc,vfc} for details. We need a  certain moduli stack, $\mathcal M_{g,E,\gamma_{v}}$,  thought of as a compactified moduli stack of genus $g$, energy $E$ holomorphic curves with $m$ marked points, where $v$ is $m$-valent.   The $m$ marked points are labeled by the $m$ possible oriented edges leaving $v$ --- so any self edge attached to $v$ labels two marked points, and all other edges attached to $v$ label a single marked point.
The curves in $\mathcal M_{g,E,\gamma_{v}}$ obey further conditions depending on $\gamma_{v}$,  described below.

\begin{itemize}
\item If $v$ is in the corner of $\triangle$ corresponding to $M_{1}$,  $\mathcal M_{g,E,\gamma_{v}}$ concerns curves in $M_{1}$.  Given any oriented edge leaving $v$ with derivative $(a,b)$, the corresponding marked point is required to intersect $M_{2}\cap M_{1}$ to degree $a$, and $M_{3}\cap M_{1}$ to degree $b$. 
\item If $v$ is in the corner corresponding to $M_{2}$, $\mathcal M_{g,E,\gamma_{v}}$ concerns curves in $M_{2}$.  Given any oriented edge leaving $v$ with derivative $(-a-b,b)$, the corresponding marked point is required to intersect $M_{1}\cap M_{2}$ to degree $a$, and $M_{2}\cap M_{3}$ to degree $b$.  
\item If $v$ is in the corner corresponding to $M_{3}$, $\mathcal M_{g,E,\gamma_{v}}$ concerns curves in $M_{3}$.  Given any oriented edge leaving $v$ with derivative $(b,-a-b)$, the corresponding marked point is required to intersect $M_{1}\cap M_{3}$ to degree $a$, and $M_{2}\cap M_{3}$ to degree $b$.  

The unified description of the above three cases is as follows: Each divisor in $M_{i}$ corresponds to an edge of $\triangle$ leaving the corner corresponding to $M_{i}$. Simple contact with this divisor corresponds to the primitive integral vector pointing away from the corner along this edge. 

\item If $v$ is on the edge corresponding to $M_{i}\cap M_{j}$, then $\mathcal M_{g,E,\gamma_{v}}$ concerns curves in (a compactification of) the $\mathbb C^{*}$--bundle $L_{ij}$ over $M_{i}\cap M_{j}$ equal to the normal bundle of $M_{i}\cap M_{j}\subset M_{i}$ minus its zero section.\footnote{The reader should not be concerned that $L_{ij}\neq L_{ji}$.} Let $\alpha$ be the primitive integral vector pointing from the corner of the dual intersection complex corresponding to $M_{i}$ to the corner corresponding to $M_{j}$, and let $\beta$ be the primitive integral vector pointing from $M_{i}$ to the remaining corner which corresponds to $M_{k}$. Then an edge leaving $v$ with derivative $a\alpha+b\beta$ corresponds to a marked point at which there is is a zero of order $a$, and this marked point is  also required to intersect the divisor corresponding to $M_{1}\cap M_{2}\cap M_{3}$ to order $b$.  When $a$ is negative, we instead require a pole of order $\abs a$.
\item If $v$ is in the interior of the dual intersection complex, then $\mathcal M_{g,E,\gamma_{v}}$ concerns curves in (a compactification of) $L_{1}\oplus L_{2}\rvert_{M_{1}\cap M_{2}\cap M_{3}}$. An edge exiting $v$ with derivative $(a,b)$ corresponds to a marked point required to have  zeros/poles in the $(L_{1},L_{2})$--direction of order $(a,b)$.
\end{itemize}
In each of the above cases, the curves in $\mathcal M_{g,E,\gamma_{v}}$ do not intersect the divisors apart from at the marked points, so their intersection with the divisors is prescribed by $\gamma$, or more specifically, the derivative of all the edges in $\gamma$ leaving $v$. Recall that we have a K\"ahler form on $M_{i}$ to define the energy of curves; the energy of curves in $L_{ij}$ and $L_{12}\oplus L_{13}\rvert_{M_{1}\cap M_{2}\cap M_{3}}$ is defined using the energy of the corresponding projections to $M_{i}\cap M_{j}$ and $M_{1}\cap M_{2}\cap M_{3}$ respectively.

\

There is an evaluation map,
\[ev:\mathcal M_{g,E,\gamma_{v}}\longrightarrow \prod_e \bar X_{e}\]
where the product is over all choices of oriented edges $e$ leaving $v$. If we label some of these $e$ outgoing, and some incoming, we get the following two evaluation maps.
\[ev_{in}:\mathcal M_{g,E,\gamma_{v}}\longrightarrow \prod_{\text{incoming }e} \bar X_{e}\]
\[ev_{out}:\mathcal M_{g,E,\gamma_{v}}\longrightarrow \prod_{\text{outgoing }e'} \bar X_{e'}\]
Then, with this designation of incoming edges $e$ and outgoing edges $e'$, $\eta_{g,E,\gamma_{v}}$ is
\[\eta_{g,E,\gamma_{v}}:=(ev_{out})_{!}\circ (ev_{in})^{*}:\rh^{*}\lrb{\prod\bar X_{e},\mathbb R}\longrightarrow \rh^{*}\lrb{\prod\bar X_{e'},\mathbb R}\]
Where $\rh^{*}$ indicates refined cohomology, discussed shortly.
Note that $\prod_{e'}\bar X_{e'}$ is the quotient of some smooth manifold $X$ by a trivial $G$--action. In the case that $\abs G$ is finite, pushforward corresponds to pushforward to $X$ multiplied by $\abs G$. In the case that $\abs G$ is infinite,  one of the internal edges of $\gamma$ has derivative $0$, and the contribution of $\gamma$ to our Gromov--Witten invariants will be $0$.
The zero-genus invariant is  defined as follows:
\[\eta_{0,\gamma_{v}}:=\sum_{E}\eta_{0,E,\gamma_{v}}q^{E}\]
and the invariant for arbitrary genus is
\[\eta_{\gamma_{v}}:=\sum_{E,g}\eta_{g,E,\gamma_{v}}q^{E}x^{2g-2+m}\]
where the valence of $v$ is $m$.

\section{Refined cohomology}\label{refined cohomology}

In our gluing formula, we are often dealing with a compactification $\bar X$ of a complex manifold $X$, where $\bar X\setminus X$ is a normal crossing divisor. For example, $X$ might be $(\mathbb C^{*})^{n}$ and $\bar X$ some compact $n$-dimensional toric manifold. If $n\geq 2$, then there are infinitely many such toric compactifications of $(\mathbb C^{*})^{n}$. The refined cohomology of $\bar X$  is the direct limit of the cohomology of all such smooth toric compactifications.  In the general case, $\bar X$ with its boundary divisor locally looks like an open subset of a toric manifold $Y$ with its toric boundary divisor. Toric blowups $Y'\longrightarrow Y$ of $Y$ correspond to subdividing the dual toric fan of $Y$. Say that a map of manifolds with normal crossing divisors $\bar X'\longrightarrow \bar X$ is a boundary blowup if it is locally modeled on toric blowups. The refined cohomology of $\bar X$ is the direct limit of the cohomology of all such boundary blowups. Refined cohomology was defined rather differently in \cite{dre}, but toric resolution of singularities, and the invariance in families of the homology theory defined in \cite{dre} implies that the two definitions are equivalent.

In practice, the boundary blowup $\bar X'$ of $\bar X$ required to faithfully  record relative Gromov--Witten invariants just needs to have the evaluation map $ev:\mathcal M\longrightarrow \bar X$ transverse to each boundary divisor. This evaluation map is actually a map in the category of exploded manifolds or log schemes; a perturbation of $ev$ within this category will not change its contact data with the divisor.

Note that in the case that the divisor in  $\bar X$ has no singularities, there are no nontrivial boundary blowups of $\bar X$.  So, in this case $\rh^{*}(\bar X)=H^{*}(\bar X)$. When the dual intersection complex has dimension at most $2$, the $\bar X_{e}$ that appear in the gluing formula have smooth boundary divisor; this is the reason that we could use ordinary cohomology for the zero-genus gluing formula in this case. An equivalent gluing formula for higher genus invariants will sometimes require a fiber product of relative moduli spaces over the product of at least two $\bar X_{e}$; the divisor in this product will no longer be smooth, and will admit infinitely many boundary blowups.

There are other cases in which ordinary cohomology can be used ---  $\rh^{0}(\bar X)=H^{0}(\bar X)$, and the same holds for top-dimensional cohomology, so when we can restrict to zero or top  dimensional cohomology, refined cohomology is not necessary; see for example the simplified gluing formula used in \cite{tec}.

\

Boundary blowups are relevant for a second reason: The relative Gromov--Witten invariants of any boundary blowup of  $\bar X$ can be identified with the relative Gromov--Witten invariants of $\bar X$ --- proved in \cite{ilgw} and \cite{egw} in the log and exploded settings respectively. This   is the reason we were not overly concerned about choosing compactifications when describing the relative Gromov--Witten invariants in section \ref{relative GW section}. 

\section{Example: $\mathbb CP^{2}$}\label{CP2}
There is a degeneration of $\mathbb CP^{2}$ into $3$ pieces, each isomorphic to $\mathbb CP^{2}$ blown up at $1$ point. In this case, $M_{1}\cap M_{2}\subset M_{1}$ is a line minus the exceptional divisor, but $M_{1}\cap M_{2}\subset M_{2}$ is the exceptional divisor, and the same statement holds for $(1,2)$ replaced by $(2,3)$ and $(3,1)$. A picture of this using toric moment maps is as follows:

\includegraphics{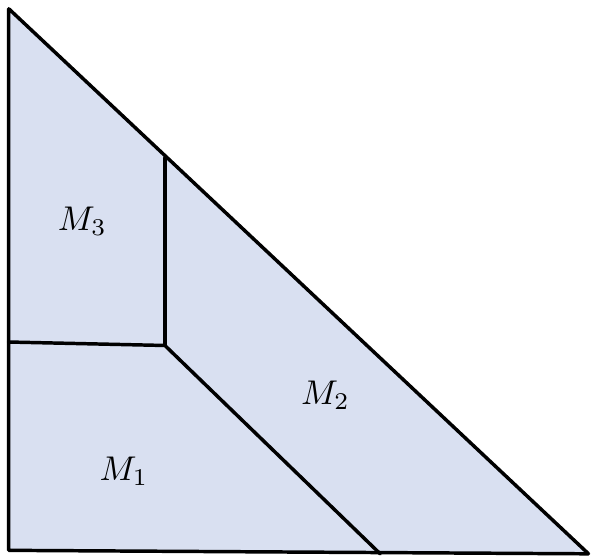}

If the points $p_{i}$ are placed in the interior of the dual intersection complex $\Delta$, then we can recover Mikhalkin's correspondence between counts of tropical curves (in $\mathbb R^{2}$ instead of the dual intersection complex) and Gromov--Witten invariants \cite{Mikhalkin}, as well as Markwig and Rau's extension of this to $0$-genus  descendant Gromov--Witten invariants, \cite{markwig}. From this perspective it should come as no surprise that tropical curves or some combinatorial analogue need  be considered to write a correct gluing formula. 

The reader might think that these tropical curves and relative invariants from vertices in the interior of $\triangle$ were only needed because we chose the points in the interior of $\triangle$, which corresponds to the un-generic seeming choice of putting all points at the intersection of all $M_{i}$. As the next example shows, tropical curves are forced upon us, even if we choose points `generically' in the interior of $M_{i}$.

 Consider degree $3$  rational curves constrained to pass through $8$ points, generically distributed $3$ in $M_{1}$, $3$ in $M_{2}$ and $2$ in $M_{3}$. In the dual intersection complex, we put these points at the corresponding corners of $\triangle$. 
The following is a picture of a tropical curve in $\triangle$  --- the little number at a  vertex in the corner denotes the number of external edges attached to that vertex, because all of these external edges are sent to a point in $\triangle$. In what follows, we shall compute the contribution of this tropical curve to Gromov--Witten invariants using our gluing formula.

\includegraphics{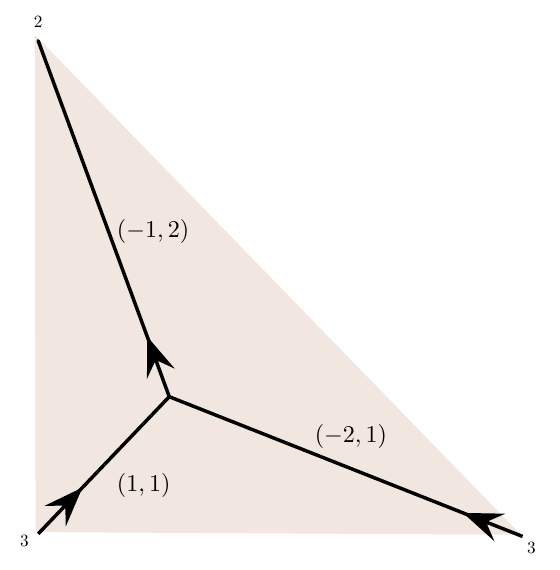}

Consider $\eta_{0,\gamma_{v}}(\tau_{0}(p_{1})\tau_{0}(p_{2})\tau_{0}(p_{3}))$ for $v$ in the bottom lefthand corner ---  this equals $\theta_{e_{1}}$ for the edge $e_{1}$ leaving that corner. We are looking for rational curves in $M_{1}$ constrained to pass through $3$ points, and with a $4$th point required to contact both divisors --- a line and the exceptional divisor --- to degree $1$. So long as our $3$ points are chosen generically, there is only one such holomorphic curve, and it has energy $E_{13}+2E_{12}$, where $E_{ij}$ is the symplectic area of $M_{i}\cap M_{j}$. In this case, $\bar X_{e_{1}}$ is $\mathbb CP^{1}$, and $\theta_{e_{1}}$ is $q^{E_{13}+2E_{12}}$ times the fundamental class in $H^{2}(\bar X_{e_{1}})$. Similarly, for $e_{2}$ the edge leaving the bottom right vertex,  $\theta_{e_{2}}$ is $q^{E_{12}+2E_{23}}$ times the fundamental class in $H^{2}(\bar X_{e_{2}})$. 

Now to compute $\theta_{e_{3}}$ for $e_{3}$ the edge leaving the central vertex $v'$. Curves in $\mathcal M_{0,\gamma_{v'}}$ are rational curves  in a toric compactification of $(\mathbb C^{*})^{2}$ that have zeros/poles at marked points of order $(-1,-1)$, $(2,-1)$, and $(-1,2)$ respectively. The moduli space of such holomorphic curves is some compactification of $(\mathbb C^{*})^{2}$, and the evaluation map to $\bar X_{e}$ is some translate of the quotient map taking the quotient by the $\mathbb C^{*}$--action with weight $(a,b)$ the derivative on the edge $e$. The `incoming' evaluation map is a degree $3$ cover of $X_{e_{1}}\times X_{e_{2}}$, because $\abs{(1,1)\wedge(-2,1)}=3$. The pullback of $\theta_{e_{1}}\wedge\theta_{e_{2}}$ to $(\mathbb C^{*})^{2}$ is $3q^{3E_{12}+2E_{23}+E_{13}}$ times the fundamental class in $H^{4}_{c}$. We push this forward to obtain $\theta_{e_{3}}\in H^{*}(\bar X_{e_{3}})$ as $3q^{3E_{12}+2E_{23}+E_{13}}$ times the fundamental class. (The energy of the curves at this vertex is $0$).

It remains to compute $\eta_{0,\gamma_{v_{0}}}(\tau_{0}(p_{7})\tau_{0}(p_{8})\theta_{3})$ for $v_{0}$ the vertex in the top corner. Our relative moduli space in consists of rational curves in $M_{3}$ constrained to pass through two points, and with a third point constrained to intersect both divisors to degree $1$. There is a 1-complex-dimensional family of such curves, all with energy $2E_{13}+E_{23}$. With these point constraints, the evaluation map to $\bar X_{e_{3}}$ has degree $1$, and $\eta_{0,\gamma_{v_{0}}}(\tau_{0}(p_{7})\tau_{0}(p_{8})\theta_{3})=3q^{3E_{12}+3E_{23}+3E_{13}}$. A degree $d$ curve has energy $d(E_{12}+E_{23}+E_{13})$, so this tropical curve contributes 3 curves towards the count of $12$ degree $3$ rational curves passing through $8$ points. The remaining $9$ curves are contributed by the tropical curves pictured below.  

\includegraphics{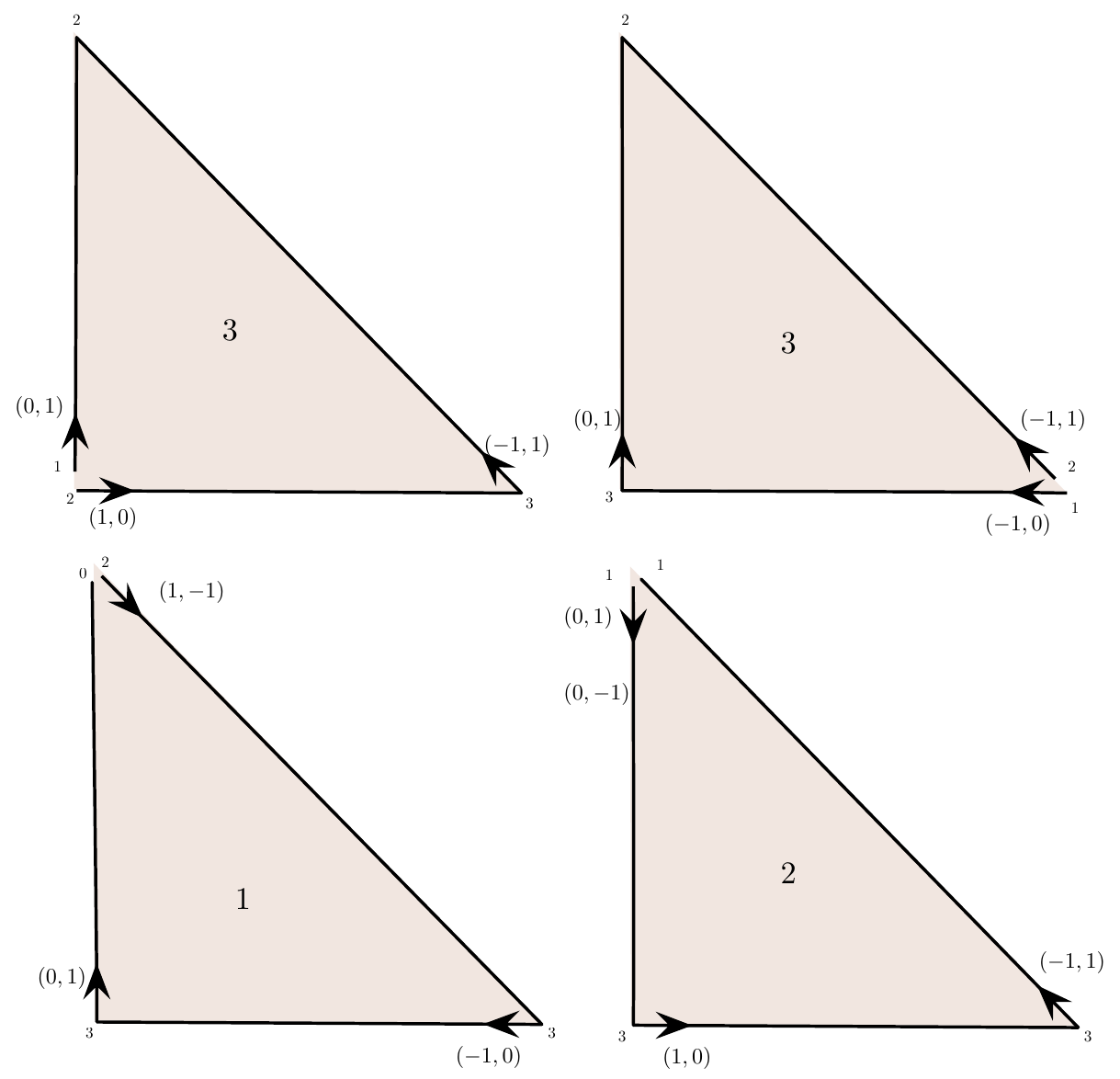}

In each of these pictures, the two vertices that are close to a corner are actually sent to that corner, but are drawn  separated to emphasize that the domain of the tropical curve is not connected there. Actually, each picture represents some  number of tropical curves corresponding to the choices of which points are attached to which of these two vertices. Each of these tropical curves contributes $q^{3E_{12}+3E_{23}+3E_{13}}$ to the count of Gromov--Witten invariants, completing the correct count of $12$ degree $3$ rational curves passing through $8$ points. 

In higher degrees, the number and types of tropical curves contributing increases rapidly. For example, in counting the degree $7$ rational curves passing through $20$ points, the following tropical curve contributes positively when $9$, $5$ and $6$ points are distributed to $M_{1}$, $M_{2}$, and $M_{3}$ respectively. From the perspective of ordinary geometry, such a  curve does not look at all generic --- it has $6$ rational components, two of which are sent entirely into the point $M_{1}\cap M_{2}\cap M_{3}$. Nevertheless, such curves are forced upon us --- they always appear in the limit of our degeneration, and contribute positively to Gromov--Witten invariants; there is no way to perturb them away, or avoid them.

\includegraphics{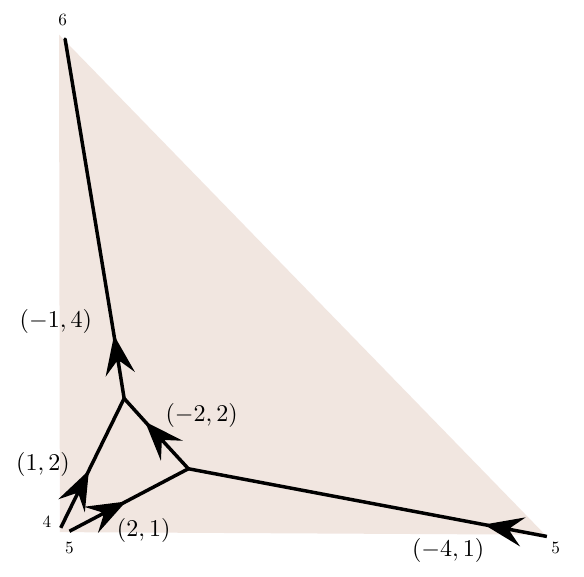}
 
\section{Example: an `unbalanced' tropical curve} \label{unbalanced}

In the previous example, all tropical curves contributing to Gromov--Witten invariants satisfy a balancing condition: for any vertex in the interior of $\triangle$,  the sum of the derivatives of edges leaving that vertex is zero. This balancing condition followed from the fact that these derivatives correspond to zeros/poles of a pair of meromorphic functions --- such a balancing condition is always  satisfied when $L_{12}\oplus L_{13}\rvert_{M_{1}\cap M_{2}\cap M_{3}}$ is a trivial bundle. In the next example, this bundle will not be trivial, and we shall see an example of an unbalanced tropical curve that contributes to Gromov--Witten invariants.  

Let us `twist' the previous example by a nontrivial line bundle as follows. The degeneration contained in the previous example can be considered within the category of toric manifolds. Let the total space of the degeneration be $X$.  Consider a $\mathbb C^{*}$--bundle $L$ over $\mathbb CP^{1}$ with Chern class $1$, and set $X'=(X\times L)/\mathbb C^{*}$ where $\mathbb C^{*}$ acts with weight $1$ on $L$ and weight $(1,-1)$ on the fibers of $X$.\footnote{We can also put a toric K\"ahler structure on $X'$ using symplectic reduction; the symplectic structure will not be relevant, so we shall ignore it.} Now $X'$ is the total space of a degeneration of some nontrivial $\mathbb CP^{2}$--bundle over $\mathbb CP^{1}$. We can make $4$ toric blowups of $X'$ so that it is a degeneration of some different toric bundle $M$ over $\mathbb CP^{1}$ whose fibers have the total moment polytope pictured below. The corresponding pieces $M_{i}$ are then bundles with fibers the smaller toric manifolds pictured below.  

\includegraphics{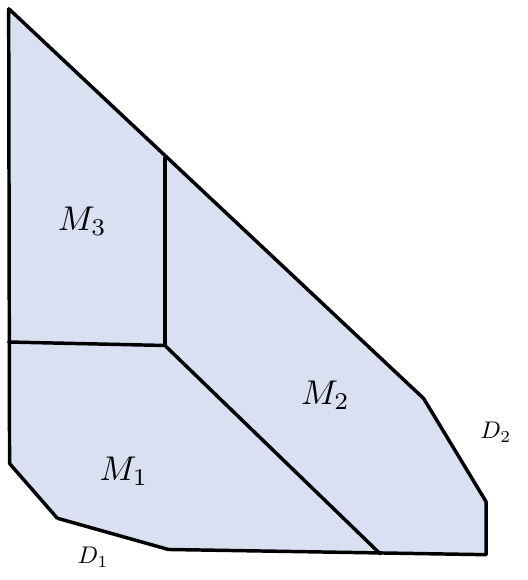}

We shall consider rational holomorphic curves in $M$ which intersect the toric boundary divisors $D_{1}$ and $D_{2}$ once,  have degree $1$ projection to $\mathbb CP^{1}$, and do not intersect any other toric boundary divisor. If  $M^{\circ}$ indicates the $(\mathbb C^{*})^{2}$--bundle over $\mathbb CP^{1}$ that is the `interior' of $M$,   we may identify a generic such curve with a meromorphic section of $M^{\circ}$ with a zero/pole of order $(-1,-2)$ and  $(2,1)$ corresponding to the intersections with $D_{1}$ and $D_{2}$ respectively. If we fix the location in $\mathbb CP^{1}$ of these singularities, the moduli space of such meromorphic sections is parametrized by $(\mathbb C^{*})^{2}$. The map from this moduli space to $D_{i}$ corresponding to the intersection with $D_{1}$ or $D_{2}$ may be thought of as recording the position of the corresponding singularity, then taking the quotient by the $\mathbb C^{*}$--action with weight $(-1,-2)$ or $(2,1)$ respectively. The corresponding map to $D_{1}\times D_{2}$ has degree $3=\abs{(-1,-2)\wedge (2,1)}$. From this we can guess that the Gromov--Witten invariant counting rational curves in this homology class constrained to pass through $2$ chosen points is $3$. 

If we choose $1$ point in $M_{1}$ and $1$ point in $M_{2}$, the only tropical curve that contributes to this Gromov--Witten invariant is pictured below.

\includegraphics{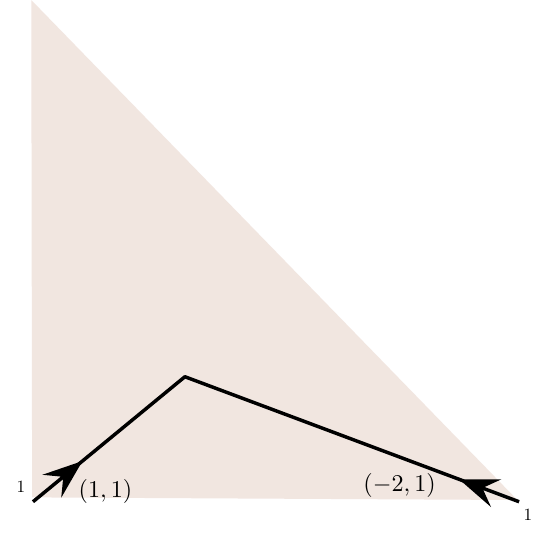}

Consider the lefthand vertex. The relative Gromov--Witten invariant at this vertex concerns rational curves in $M_{1}$ with one marked point  constrained to contact $M_{1}\cap M_{2}$ and $M_{1}\cap M_{3}$ to degree $1$, and a second marked point sent to a chosen point in $M_{1}$. The curves relevant to our computation also intersect $D_{1}$ once, do not intersect any other toric boundary divisor of $M_{1}$, and project to have degree $0$ in $\mathbb CP^{1}$. There is exactly one such curve.   Similarly, when considering the righthand vertex,  there is exactly one curve in $M_{2}$ relevant to our calculation --- it intersects $D_{2}$ once, projects to a point in $CP^{1}$, and does not intersect any other toric boundary divisor of $M_{2}$ (apart from $M_{2}\cap M_{3}$ and $M_{2}\cap M_{1}$ at a marked point).

 How about the central vertex? The relevant curves are rational curves in $L_{12}\oplus L_{13}\rvert_{M_{1}\cap M_{2}\cap M_{3}}$\footnote{In this case, we can identify $M^{\circ}$ with $L_{12}\oplus L_{13}\rvert_{M_{1}\cap M_{2}\cap M_{3}}$, however  the description of $M^{\circ}$ as a $(\mathbb C^{*})^{2}$--bundle above uses $(L_{12}^{-1}L_{13})\oplus L_{13}^{-1}\rvert_{M_{1}\cap M_{2} \cap M_{3}}$.}, with one puncture with  a zero/pole of order $(-1,-1)$, and another puncture with a zero/pole of order $(2,-1)$
. This is only possible for curves that project to $\mathbb CP^{1}=M_{1}\cap M_{2}\cap M_{3}$ to have degree $1$, so we can identify such curves with meromorphic sections of $L_{12}\oplus L_{13}\rvert_{M_{1}\cap M_{2}\cap M_{3}}$ with the given singularities. If the location of the singularities in $\mathbb CP^{1}$ is fixed, then the moduli space of such sections is equal to $(\mathbb C^{*})^{2}$. In the above tropical curve, let $e_{1}$ be the left edge and $e_{2}$ be the right edge. $X_{e_{1}}$ and $X_{e_{2}}$ is given by taking the quotient of $L_{12}\oplus L_{13}\rvert_{M_{1}\cap M_{2}\cap M_{3}}$ by the $\mathbb C^{*}$--action of weight $(1,1)$ and $(-2,1)$ respectively. The evaluation map to $X_{e_{i}}$ is given by taking note of the position of the relevant zero/pole, and then taking the quotient by the relevant $\mathbb C^{*}$--action. It follows that the evaluation map from our moduli space to $X_{e_{1}}\times X_{e_{2}}$ has degree $3=\abs{(1,1)\wedge (-2,1)}$. From (slightly more rigorous versions of) these calculations, we can conclude that this tropical curve contributes a count of $3$ curves to our Gromov--Witten invariant. Although we were not keeping track of energy in this example, note that unlike in previous examples,  these curves have positive energy at the interior vertex.
 
 \section{Why learn about exploded manifolds?}\label{exploded}
 
\

The gluing formula described above, and more powerful generalizations are naturally proved and stated using the category of exploded manifolds, introduced in \cite{iec}. Analysis on exploded manifolds is much the same as analysis on smooth manifolds, but for studying Gromov--Witten invariants, exploded manifolds have significant advantages:
\begin{itemize}
\item Normal crossing degenerations occur in smooth families within the category of exploded manifolds. In particular,
the normal crossing degeneration of $M$ into $M_{i}$ corresponds to a smooth family of exploded manifolds.
\item Gromov--Witten invariants are defined within the category of exploded manifolds, and do not change in connected smooth families. In particular, we can compute the Gromov--Witten invariants of $M$ using calculations in an exploded manifold $\bold M$ corresponding to our singular fiber. The gluing formula is such a calculation, and with an understanding of exploded manifolds, it is obvious.\footnote{Even though the gluing formula is `obvious', it takes quite a lot of work to prove, simply because it takes a lot of work to define and work rigorously with Gromov--Witten invariants.}
\item A holomorphic curve in an exploded manifold $\ex M$ is a map in the category of exploded manifolds. Its domain $\ex C$ is an exploded manifold, the target $\ex M$ is an exploded manifold, and the map  itself $f:\ex C\longrightarrow \ex M$ is a morphism in the category of exploded manifolds. This is more natural than the usual setup for relative invariants, where holomorphic buildings or rubber components appear. 
\item As well as individual holomorphic curves being described as morphisms in the category of exploded manifolds, families of curves also admit a natural definition in the category of exploded manifolds. A  smooth family of curves in $\ex M$ is a map from the total space of a smooth family of curves (in the category of exploded manifolds) to $\ex M$. Exploded manifolds have an advantage over smooth manifolds here: because normal crossing degenerations happen within smooth families of exploded manifolds, smooth families of curves in the category of exploded manifolds naturally allow for bubbling and node formation.  
 \item After the usual caveats about transversality, the moduli stack of holomorphic curves in an exploded manifold naturally has the structure of a smooth\footnote{Actually, in the almost complex case, the moduli stack has a certain level of regularity $\C\infty1$, which is as good as smooth for all practical purposes.} 
 orbifold in the category of exploded manifolds. The same holds for the moduli stack of holomorphic curves in a family of exploded manifolds. In particular, there is no need to separately consider `smoothability' of curves in an exploded manifold corresponding to singular fiber of a normal crossing degeneration.
 \item There is a tropical part functor which sends an exploded manifold to some union of integral-affine polytopes glued along faces. In our example, the tropical part of $\ex M$ is the dual intersection complex $\triangle$, and the tropical part of a holomorphic curve in $\ex M$ is a tropical curve in $\triangle$. The tropical part of the moduli stack of holomorphic curves in $\ex M$ is a kind of  cover of the moduli stack of tropical curves in $\triangle$, so consideration of tropical curves can often determine a lot of information about the moduli stack of curves --- for examples, see \cite{tec, 3d}. 
 \end{itemize}

There are now three roughly equivalent approaches to defining Gromov--Witten invariants relative normal crossing divisors. In \cite{tropicalIonel}, I sketch the relationship between Ionel's approach to GW invariants relative normal crossing divisors, \cite{IonelGW}, and my approach using exploded manifolds. In \cite{elc}, I explain why I believe that in the algebraic case, using exploded manifolds should give the same invariants as log Gromov--Witten invariants defined by Gross and Siebert in \cite{GSlogGW} and Abramovich and Chen in \cite{Chen,acgw}. In each case, only minor modifications are needed to define invariants using some equivalent of refined cohomology for use in the gluing formula. 

\

The reader interested in an introduction to exploded manifolds should read \cite{iec,scgp}, or \cite{elc} if  already conversant with log geometry. Slides from a talk with lots of pictures explaining how tropical curves appear from an analytic perspective can be found by googling `Brett Parker tropical gluing prezi'. A version of Gromov compactness for holomorphic curves in exploded manifolds is proved in \cite{cem}, and regularity of the $\dbar$ equation in families is studied in \cite{reg}, the construction of Gromov--Witten invariants is completed in \cite{uts,evc,vfc}, and the tropical gluing formula for Gromov--Witten invariants is proved in \cite{gfgw}.

\bibliographystyle{plain}
\bibliography{ref}
\end{document}